\documentclass{article}

\usepackage{url}

\usepackage{amsmath, amssymb, amstext, amsopn, amsxtra, graphicx, color, calligra, hyperref, url, ulem, polynom, permute, ifthen, float, epstopdf, linsys, marvosym, rotating}

\newtheorem{theorem}{Theorem}
\newtheorem{corollary}{Corollary}
\newtheorem{lemma}{Lemma}

\begin{document}

\newcommand{\recip}[1]{\ensuremath{\frac1{#1}}}
\newcommand{\half}{\recip2}
\newcommand{\third}{\recip{3}}
\newcommand{\twothirds}{\ensuremath{\frac23}}
\newcommand{\fourth}{\recip4}
\newcommand{\fifth}{\recip5}
\newcommand{\sixth}{\recip6}
\newcommand{\point}[1]{\mbox{\(\left(#1\right)\)}}
\newcommand{\qar}{\begin{eqnarray*}}
\newcommand{\raq}{\end{eqnarray*}}
\newcommand{\inv}{^{-1}}
\newcommand{\vone}{{\bf 1}}
\newcommand{\va}{{\bf a}}
\newcommand{\vb}{{\bf b}}
\newcommand{\vr}{{\bf r}}
\newcommand{\vs}{{\bf s}}
\newcommand{\vu}{{\bf u}}
\newcommand{\comment}[1]{\hskip 5mm \parbox{3in}{#1}}
\newcommand{\casI}{\renewcommand{\labelenumi}{\bf Case \Roman{enumi}.}\num}
\newcommand{\sacI}{\mun \renewcommand{\labelenumi}{\arabic{enumi}.}}
\newcommand{\num}{\begin{enumerate}}
\newcommand{\mun}{\end{enumerate}}
\newcommand{\rowr}{\ensuremath{\left(\begin{array}{rrrr}r_1 & r_2 & \cdots & r_k \end{array}\right)}}
\newcommand{\rows}{\ensuremath{\left(\begin{array}{rrrr}s_1 & s_2 & \cdots & s_k \end{array}\right)}}
\newcommand{\qed}{\mbox{\(\square\,\,\)}\\}
\newcommand{\longcomment}[1]{}

\newcommand{\nakaut}{Nakayama automorphism}
\newcommand{\nakauts}{Nakayama automorphisms}
\newcommand{\scg}{\ensuremath{\dk/\dk^2}}
\newcommand{\Rd}{\ensuremath{\hat R}}
\newcommand{\Hom}[1]{\mbox{\rm Hom}_{#1}}
\newcommand{\leftmod}{\mbox{left $R$-module}}
\newcommand{\leftmods}{\mbox{left $R$-modules}}
\newcommand{\rightmod}{\mbox{right $R$-module}}
\newcommand{\rightmods}{\mbox{right $R$-modules}}
\newcommand{\Rbar}{\ensuremath{\bar R}}
\newcommand{\R}{{_R}}
\newcommand{\lsoc}{\mbox{\rm soc}(\R R)}
\newcommand{\rsoc}{\mbox{\rm soc}(R_R)}
\newcommand{\lcat}{\ensuremath{\R\M}}
\newcommand{\rcat}{\ensuremath{\M_R}}
\newcommand{\bimod}{$(R,R)$-bimodule}
\newcommand{\bimods}{$(R,R)$-bimodules}
\newcommand{\Aut}{\mbox{\rm Aut}}
\newcommand{\End}{\mbox{\rm End}}
\newcommand{\dk}{\ensuremath{\dot{k}}}
\newcommand{\M}{\mathfrak M}
\newcommand{\m}{\ensuremath{\mathfrak m}}
\newcommand{\isom}{\simeq}
\newcommand{\isommap}{\stackrel{\sim}{\longrightarrow}}
\newcommand{\lomar}{\cite{lomar}}
\newcommand{\Rrad}{\ensuremath{R/\mbox{rad}R}}
\newcommand{\rRd}{\ensuremath{\hat R_R}}
\newcommand{\lRd}{\ensuremath{\R\hat R}}
\newcommand{\lRda}{\ensuremath{\R\hat R_1}}
\newcommand{\lRdb}{\ensuremath{\R\hat R_2}}
\newcommand{\rRda}{\ensuremath{(\hat R_1)_R}}
\newcommand{\rRdb}{\ensuremath{(\hat R_2)_R}}
\newcommand{\lRdi}{\ensuremath{\R\hat R_i}}
\newcommand{\rRdi}{\ensuremath{(\hat R_i)_R}}
\newcommand{\Xd}{\ensuremath{\hat{X}}}
\newcommand{\rXd}{\ensuremath{\hat{X}_R}}
\newcommand{\rXdi}{\ensuremath{(\hat{X}_i)_R}}
\newcommand{\biRd}{\ensuremath{\R\hat R_R}}
\newcommand{\biRdi}{\ensuremath{\R(\hat R_i)_R}}
\newcommand{\biRda}{\ensuremath{\R(\hat R_1)_R}}
\newcommand{\biRdb}{\ensuremath{\R(\hat R_2)_R}}
\newcommand{\hi}{\ensuremath{\varphi}}
\newcommand{\sig}{\ensuremath{\sigma}}
\newcommand{\alp}{\ensuremath{\alpha}}
\newcommand{\el}{\ensuremath{\ell}}
\newcommand{\ro}{\ensuremath{\rho}}
\newcommand{\rad}{\ensuremath{\mbox{rad }R}}
\newcommand{\all}{ \hskip 5mm (\forall}
\newcommand{\Id}{\mbox{\rm Id}}
\newcommand{\Tr}{\mbox{\rm Tr}}
\newcommand{\field}[1]{\ensuremath{\mathbb #1}}
\newcommand{\CC}{\field C}
\newcommand{\F}{\mathcal F}
\newcommand{\G}{\mathcal G}

\newtheorem{thm}{Theorem}
\newtheorem{cor}[thm]{Corollary}
\newtheorem{conj}[thm]{Conjecture}
\newtheorem{lem}[thm]{Lemma}
\newtheorem{rem}[thm]{Remark}
\newtheorem{exmp}[thm]{Example}

\renewcommand{\em}{\textit}

\title{Nakayama automorphisms of Frobenius algebras}

\author{Will Murray\\
wmurray@csulb.edu\\
Department of Mathematics\\
California State University, Long Beach\\
Long Beach, CA 98040-1001}
\date{February 24, 2014}

\maketitle

\begin{abstract}
We show that the Nakayama automorphism of a Frobenius algebra $R$ over a field $k$ is independent of the field (Theorem~\ref{theorem-nakind}).  Consequently, the $k$-dual functor on \leftmods\ and the bimodule isomorphism type of the $k$-dual of $R$, and hence the question of whether $R$ is a symmetric $k$-algebra, are independent of $k$.  We give a purely ring-theoretic condition that is necessary and sufficient for a finite-dimensional algebra over an infinite field to be a symmetric algebra (Theorem~\ref{theorem-symring}).  

\end{abstract}

\section{Introduction}

Let $R$ be a finite-dimensional algebra over a field $k$.  The $k$-dual $\Rd := \Hom{k}(R,k)$ has a natural structure as an \bimod.  We say $R$ is a {\em Frobenius algebra} if $R \isom \Rd$ as \leftmods, and $R$ is a {\em symmetric $k$-algebra} if $R \isom \Rd$ as \bimods.  It is well-known that $\Rd$ is isomorphic to the injective hull of $\Rrad$ as \leftmods, so $R$ is Frobenius iff $R \isom E(\R(\Rrad))$ as \leftmods.  This purely ring-theoretic criterion shows that the property of $R$ being Frobenius is independent of the field $k$ over which we are considering $R$ as an algebra.  Motivated by this property,  an arbitrary artinian ring $S$ is defined to be a {\em Frobenius ring} if \(S \isom E(_S(S/\mbox{rad }S))\) as left $S$-modules, and this definition has led to a rich theory of Frobenius rings (see, for example, Section 16 in \lomar) that is not dependent on the framework of linear algebra.  

The facts above naturally raise several questions.  Is the property of $R$ being a symmetric $k$-algebra independent of $k$?  If $R$ is symmetric, we know by Brauer's Equivalence Theorem (16.70 in \lomar) that the $k$-dual functor $\Hom{k}(-,k)$ from \leftmods\ to \rightmods\ is independent of $k$, i.e. the two functors defined by different fields are naturally equivalent.  On the level of modules, this means that the \rightmod\ isomorphism type of the $k$-dual of any left module $\R X$ is independent of $k$.  Do these facts remain true if $R$ is only Frobenius?  The result above shows only that the isomorphism type of the dual of the left regular module $\R R$ is independent of $k$.

The key to all of these questions is the Nakayama automorphism, a distinguished $k$-algebra automorphism of a Frobenius algebra $R$ that measures how far $R$ is from being a symmetric algebra.  (The automorphism is the identity iff $R$ is symmetric.)  We will show that the Nakayama automorphism is independent of $k$ and derive affirmative answers to the questions above as corollaries.  We will give a purely ring-theoretic condition that is equivalent to the property of $R$ being symmetric at least in the case when $k$ is infinite.  We hope that this will promote a ring-theoretic development of properties of symmetric algebras that parallels the theory of Frobenius rings.  

F. G. Frobenius himself pioneered the idea of comparing an algebra with its dual in \cite{frobenius}.  The main properties of Frobenius algebras and symmetric algebras were developed by Nakayama in \cite{nak1}, \cite{nak2}, and \cite{nak3}.  They have been the subject of continued interest because of connections to such diverse areas as group representations, topological quantum field theories, Gorenstein rings in commutative algebra, Hopf algebras, coding theory, and the Yang-Baxter Equation.  For an excellent reference on the subject, see \lomar.

\section{The Nakayama automorphism}

In this section we show that the \nakaut\ of a Frobenius algebra is independent of the ground field.  As a corollary to the proof, we derive a simple ring-theoretic characterization of local symmetric algebras.  

Let $R$ be a finite-dimensional algebra over a field $k$.  In \lomar, Theorem 3.15, we have:

\begin{theorem} \label{Frobdef}The following are equivalent:
\begin{description}
\item[1.]  $R$ is a Frobenius algebra, i.e.\  \(R \isom \Rd\) as \leftmods.
\item[2.]  There exists a linear functional \(\lambda: R \rightarrow k\) whose kernel contains no nonzero 
left ideals.
\item[3.]  There exists a hyperplane $H \subset R$ (i.e. a subspace of codimension 1) containing no nonzero left ideals.
\item[4.]  There exists a nondegenerate associative bilinear form \(B: R \times R \rightarrow k\).  (``Associative'' means \(B(rs,t) = B(r,st)\).)
\end{description}
\end{theorem}  

The equivalence of the first two conditions follows from taking \(\lambda\) to be the image of 1 under the module isomorphism and vice versa.  The equivalence of the second and fourth condition follows from defining \(B(r,s):=\lambda(rs)\) and \(\lambda(r):= B(r,1)\).  Since the last condition is right-left symmetric, we could also include the right-handed analogues of the other conditions above.  

Given one isomorphism \(\hi:R \isommap \Rd\), any other isomorphism $\hi'$ is obtained by composition with an automorphism of the left regular module $_R R$, which corresponds to right multiplication by a unit $u \in U(R)$.  This affects the other conditions above as follows:  the new functional is \(\lambda' = u\lambda: r \mapsto \lambda(ru)\); the new hyperplane is \(H' = \ker \lambda' = Hu\inv\); and the new form is \(B'(r,s) = B(r,su)\).

A similar theorem (\lomar, Theorem 16.54) applies to symmetric $k$-algebras:

\begin{theorem} \label{symdef}The following are equivalent:
\begin{description}
\item[1.]  $R$ is a symmetric algebra, i.e.\  \(R \isom \Rd\) as \bimods.
\item[2.]  There exists a functional \(\lambda: R \rightarrow k\) such that \(\ker \lambda\) contains no nonzero left ideals and \(\lambda(rs) = \lambda(sr) \; \forall r,s \in R\).
\item[3.]  There exists a hyperplane $H \subset R$ containing the commutators \([R,R] = \{\sum_i (r_is_i - s_ir_i):r_i, s_i \in R\}\) and containing no nonzero left ideals.
\item[4.]  There exists a nondegenerate associative symmetric bilinear form \(B: R \times R \rightarrow k\).
\end{description}
\end{theorem}  

If the conditions of Theorem~\ref{Frobdef} hold, the nondegeneracy of the form $B$ implies that there is a unique $k$-linear map \(\sigma: R \rightarrow R\) defined by \(B(r,s)=B(s,\sigma(r)) \; \forall r,s \in R\).  It is easy to check that $\sigma$ is actually a $k$-algebra automorphism of $R$; we call it the {\em \nakaut} of $R$.  Replacing $B$ with a new form $B'$ defined by the unit $u$ gives us the new automorphism \(\sigma':r \mapsto u\sigma(r)u\inv\).  So the \nakaut\ is determined up to composition with inner automorphisms; equivalently, it is a well-defined element of the group of outer automorphisms of $R$.  The algebra is symmetric iff \sig\ can be taken to be the identity, iff the \nakaut\ determined by an arbitrary nondegenerate associative bilinear form is an inner automorphism.

If we use the linear functional $\lambda$ to define \sig\ instead of the form $B$, then \sig\ is defined by the equation

\[
\lambda(rs) = \lambda(s\sigma(r)) \all r,s \in R).
\]

We are now ready to prove that the \nakaut\ is independent of the base field.  We warm up with the local case.  The argument is similar to that for the general case but much easier, and it gives us a criterion for a local algebra to be symmetric.  

\begin{theorem}\label{theorem-local}
If $R$ is a local Frobenius k-algebra then $\sigma$ is independent of $k$.
\end{theorem}

\textit{Proof}.
Let $k_1$ and $k_2$ be two fields over which $R$ is a finite-dimensional algebra, and suppose $\sigma_1$ is a Nakayama automorphism of $R$ as a $k_1$-algebra.  Then $\sigma_1$ arises from a $k_1$-linear functional \(\lambda_1:R \rightarrow k_1\) via the equation
\[
\lambda_1(rs) = \lambda_1(s \sigma_1(r)) \all r,s \in R).  
\]

Thus \(C:= \{\sum(r_is_i - s_i\sigma_1(r_i)) :\; r_i, s_i \in R\} \subseteq \ker \lambda_1\).  Note that $C$ is closed under multiplication by any element from the center $Z(R)$, and in particular that $C$ is a subspace with respect to both $k_1$ and $k_2$.  

Now since $R$ is local Frobenius, $\R R$ is the only principal indecomposable \leftmod, and so $\R R$ has a simple socle $S$ by Theorem 16.4 in \lomar.  Then \(S \not\subset \ker \lambda_1\), so 
\(
S \not\subset C.  
\)

Since $S$ and $C$ are both $k_2$-subspaces, we can define a $k_2$-linear functional \(\lambda_2:R \rightarrow k_2\) that is $0$ on $C$ but not on $S$.  Then since \(S \not\subset \ker\lambda_2\), \(\ker\lambda_2\) contains no nonzero left ideals, and the Nakayama automorphism $\sigma_2$ of $R$ as a Frobenius $k_2$-algebra is defined by 
\[
\lambda_2(rs) = \lambda_2(s \sigma_2(r)) \all r,s \in R).  
\]
In other words, $\sigma_2(r)$ is uniquely defined by 
\[
rs - s \sigma_2(r) \in \ker\lambda_2 \all s \in R).  
\]
But \( rs - s \sigma_1(r) \in C \subseteq \mbox{ ker}\lambda_2 \; (\forall s)\), so \(\sigma_2(r) = \sigma_1(r) \; \forall r \in R\), as desired.
  \qed

The proof above gives us the promised ring-theoretic characterization of local symmetric algebras.  Recall that the property of $R$ being Frobenius over $k$ is independent of $k$, and in fact, is equivalent to a ring-theoretic property.

\begin{corollary}\label{corollary-localsym}
Let $R$ be a local $k$-algebra.  Then $R$ is a symmetric $k$-algebra iff $R$ is a Frobenius $k$-algebra and \( \lsoc \not\subset [R,R]\).  In particular, the truth of $R$ being a symmetric $k$-algebra is independent of $k$.
\end{corollary}

\textit{Proof}.
This follows from the proof of the theorem above.  If $R$ is symmetric, then we can take $\sigma_1$ to be the identity, so \(S \not\subset C = [R,R]\).  Conversely, if \(S \not\subset [R,R]\), then we can define $\lambda_2$ as we did above to be $0$ on $[R,R]$ but not on $S$.  The resulting $\sigma_2$ will be the identity, proving that $R$ is a symmetric algebra.
  \qed

We now pass to the general case and show that the \nakaut\ with respect to the two fields remains the same.  This turns out to be easy if the fields are both finite-dimensional over their intersection (necessarily a field).   The case in which there is no convenient intersection is harder and uses the assumption that the fields be infinite, so we do not have a single proof to cover both cases.

Let $R$ be a Frobenius ring with Jacobson radical $J$ and \(\Rbar = R/J\).  Suppose, as above, that $R$ can be considered as a finite-dimensional algebra over two different fields $k_1$ and $k_2$, with respective \nakauts\ $\sigma_1$ and $\sigma_2$.  

\begin{theorem}  \label{theorem-nakind}The Nakayama automorphism of $R$ is independent of the ground field.
\end{theorem}

\textit{Proof of Part I}.
Assume that the two fields are both finite-dimensional over some common ground field.  This case can be handled by a transfer-type argument, as suggested to me by T.Y. Lam.  By passing down to the common ground field and then up again, we can reduce to the case in which \(k_2 \subseteq k_1\).  

Let \(\mbox{Tr}:k_1 \rightarrow k_2\) be any nonzero $k_2$-linear map.  Considering $R$ as a Frobenius $k_1$-algebra, we have a $k_1$-linear functional \(\lambda_1:R \rightarrow k_1\) whose kernel contains no nonzero left ideals.  Then \(\lambda_2 := \Tr \circ \lambda_1: R \rightarrow k_2\) is a $k_2$-linear map, and we claim that \(\ker \lambda_2\) also contains no nonzero left ideals.  Indeed, if $r \in R\setminus\{0\}$, then \(\exists s \in R\) such that \(\lambda_1 (sr) \neq 0\), so \(\exists \alpha \in k_1\) such that \(0 \neq \Tr (\alpha \lambda_1(sr)) = \Tr (\lambda_1 (\alpha sr)) = \lambda_2 (\alpha sr).\)  

Now $\sigma_i(r)$ is defined (\(\forall r \in R\)) by the equation
\[
rs - s \sigma_i(r) \in \ker \lambda_i \all s \in R).  
\]
But \(\ker \lambda_1 \subseteq \ker \lambda_2\), so (\(\forall r \in R\)), 
\[
rs - s \sigma_1(r) \in \ker \lambda_1 \subseteq \ker \lambda_2 \all s \in R), 
\]
showing that $\sigma_2(r)$ must be equal to $\sigma_1(r)$.  This finishes Part I\@.
  \qed

For Part II we first need two facts from linear algebra.

\begin{lemma}\label{hyperplane}
Let $U \subsetneq V$ be finite-dimensional vector spaces over a field $k$ and suppose $V$ decomposes into subspaces \(V = V_1 \oplus V_2 \oplus \cdots \oplus V_n\) with each \(V_i \not\subset U\).  Suppose that $|k| \geq n$.  Then $U$ can be enlarged to a hyperplane $U'$ such that \(V_i \not\subset U'\) for \(i = 1,2,\dots, n\).
\end{lemma}

\textit{Proof}.
By enlarging $U$ one dimension at a time, we may assume that $U$ is maximal with respect to the property that no \(V_i \subseteq U\).  We claim that now \(\dim_k V/U = 1\).  If not, there exist 
at least $|k| + 1$ linear subspaces of $V/U$ corresponding to one-dimensional extensions $U_i \supset U$.  By the maximality of $U$ and the Pigeonhole Principle, some $V_i$ is contained in two different extensions, say $U_1$ and $U_2$.  But this implies \(V_i \subseteq U_1 \cap U_2 = U\), a contradiction.
  \qed

(The assumption that $|k| \geq n$ cannot be omitted.  A three-dimensional vector space over the field of two elements contains the subspace \(U = \{0, (1,1,1)\}\), which cannot be extended to a hyperplane without including one of the three coordinate axes.)

\begin{lemma}\label{lemma-commutators}
Let $D$ be a division ring, $n$ a positive integer, and $S = \mathbb M_n(D)$.  If $I \subseteq S$ is any nonzero left ideal, then $I + [S,S] = S$.  
\end{lemma}

\textit{Proof}.
Let $U = I + [S,S]$ and let $E_{ij}$ denote the matrix units in $S$.  Using a nonzero element of $I$, we can obtain a matrix in $U$ that is nonzero in the $(i,i)$ position and $0$ off the $i$-th row.  For all $d \in D$ and $i \neq j$, \(dE_{ij} =  (dE_{ii})(E_{ij}) - (E_{ij})(dE_{ii}) \in U\), and \(d(E_{ii} - E_{jj}) = (dE_{ij})(E_{ji}) - (E_{ji})(dE_{ij}) \in U\).  Repeated use of these identities shows that an arbitrary matrix in $S$ is a sum of matrices in $U$.
  \qed

\textit{Proof of Theorem~\ref{theorem-nakind}, Part II}.  We assume that there is no common ground field over which $k_1$ and $k_2$ are both finite-dimensional.  We need this assumption only because we will need to assume that both fields are infinite so that we can apply Lemma~\ref{hyperplane}.  

Fix a $k_1$-linear functional \(\lambda_1: R \rightarrow k_1\) with kernel $H_1$ containing no nonzero left ideals.  Then the \nakaut\ of $R$ as a $k_1$-algebra is defined (\(\forall r \in R\)) by 
\[
rs - s \sigma_1(r) \in H_1 \all s \in R).
\]
As in the proof of Theorem~\ref{theorem-local}, we set
\[
C:= \left\{\sum_i(r_is_i - s_i\sigma_1(r_i)) :\; r_i, s_i \in R\right\} \subseteq H_1
\]
and note that $C$ is closed under multiplication from the center $Z(R)$.  In particular, $C$ is a subspace over both $k_1$ and $k_2$.  Let \(S:= \lsoc\) and note that since \(C \subseteq H_1\), $S \cap C$ contains no nonzero left ideals.  Also \(S \cap C \subseteq S \cap H_1\), which is a $k_1$-subspace of $S$ of codimension 1.  (This is because \(\dim_{k_1} R/H_1 = 1\) and \(S \cap H_1 \neq S\) because $H_1$ contains no nonzero left ideals.)  

By Theorem 16.14 in \lomar, we have an isomorphism \(\hi: \R S \isommap \R \Rbar\), which is also an isomorphism of left \Rbar-modules.  Now \(S \cap H_1 \subset S\) contains no nonzero left ideals of $R$, hence no minimal left ideals, hence no nonzero \Rbar-submodules.  So \(\hi(S \cap H_1)\) is a $k_1$-hyperplane in \Rbar\ containing no nonzero left ideals.  

Since $R$ is a finite-dimensional algebra (over either field), \Rbar\ is semisimple (by Theorem 4.14 in \cite{fc}), hence a symmetric algebra by Example 16.59 in \lomar.  We consider \Rbar\ now as a symmetric $k_1$-algebra.  By Theorem~\ref{symdef}, \Rbar\ contains another $k_1$-hyperplane $H$ that contains no nonzero left ideals and contains the commutator subspace \([\Rbar, \Rbar]\).  Now by the discussion following Theorem~\ref{Frobdef}, we know that \(H = (\hi(S \cap H_1))u\) for some \(u \in U(\Rbar).\) 

Now \((\hi(S \cap C))u \subseteq \hi(S \cap H_1))u = H\), so \(U:=(\hi(S \cap C))u + [\Rbar, \Rbar] \subseteq H\).  Since $H$ contains no nonzero left ideals in \Rbar, $U$ also contains no nonzero left ideals.  But since both \(\hi(S \cap C))u\) and \([\Rbar, \Rbar]\) are $k_2$-subspaces of $\Rbar$, $U$ is a $k_2$-subspace of \Rbar.  Our goal is to enlarge $U$ to a $k_2$-hyperplane containing no nonzero left ideals.  

Let \Rbar\ have Artin-Wedderburn decomposition
\(
\mathbb M_{n_1} (D_1) \times \cdots \times \mathbb 
M_{n_r} (D_r),
\)
where the $D_i$'s are division rings.  We decompose each \(R_i:= \mathbb M_{n_i} (D_i)\) into a sum of simple left ideals $V_{i,j}$, where $V_{i,j}$ consists of matrices that are 0 except in the $j$-th column.  This gives a decomposition of \Rbar\ into simple left ideals:
\[
\Rbar = V_{1,1} \oplus \cdots \oplus V_{1,n_1} \oplus \cdots \oplus V_{r,1} \oplus \cdots \oplus V_{r,n_r}.  
\]
Now we know that for all $i,j$, \(V_{i,j} \not\subset U\).  So by Lemma~\ref{hyperplane} we can enlarge $U$ to a $k_2$-hyperplane \(U' \subset \Rbar\) while preserving \(V_{i,j} \not\subset U' \forall i,j\).  We claim that $U'$ still contains no nonzero left ideal of \Rbar.  Indeed, assume that $U'$ does contain a nonzero left ideal of \Rbar; then it contains a minimal left ideal of one of the $R_i$'s, say $R_1$.  But $U'$ also contains the commutators $[R_1,R_1]$ since \([R_1,R_1] \subseteq [\Rbar, \Rbar] \subseteq U \subseteq U'\).  Then by Lemma~\ref{lemma-commutators}, $U'$ contains all of $R_1$, hence all the $V_{1,j}$'s, a contradiction.  So $U'$ is indeed a $k_2$-hyperplane of \Rbar\ containing no nonzero left ideals.  

We now consider the $k_2$-hyperplane \(U'u\inv \subset \Rbar\), which also contains no nonzero left ideals of \Rbar.  Moreover, since \((\hi(S \cap C))u \subseteq U \subseteq U'\), we have \(\hi(S \cap C) \subseteq U'u\inv\).  We now pull \(U'u\inv\) back through the isomorphism \(\hi: \R S \isommap \R \Rbar\) to get a $k_2$-hyperplane \(H_2':= \hi\inv(U'u\inv) \subset S\) containing no nonzero left \Rbar-submodules of $S$, hence no nonzero left ideals of $R$.  Also, since \(\hi(S \cap C) \subseteq U'u\inv\), \(H_2'\) contains $S \cap C$.  

To finish the proof, we will extend $H_2'$ to a $k_2$-hyperplane $H_2 \subset R$ that contains $C$ and still contains no nonzero left ideals.  We can then use $H_2$ to define the \nakaut\ with respect to $k_2$. 

\setlength{\unitlength}{.1in}

\begin{figure}
\begin{center}
\begin{picture}(23,18)(0,0)

\put(0,0){\framebox(23,18)[tl]{}}

\put(5,5){\oval(6,6)[bl]}
\put(6,5){\oval(6,6)[br]}
\put(6,13){\oval(6,6)[tr]}
\put(5,13){\oval(6,6)[tl]}
\put(5,2){\line(1,0){1}}
\put(5,16){\line(1,0){1}}
\put(2,5){\line(0,1){8}}
\put(9,5){\line(0,1){8}}

\put(5,5){\oval(4,4)[bl]}
\put(7,5){\oval(4,4)[br]}
\put(7,13){\oval(4,4)[tr]}
\put(5,13){\oval(4,4)[tl]}
\put(5,3){\line(1,0){2}}
\put(5,15){\line(1,0){2}}
\put(3,5){\line(0,1){8}}
\put(9,5){\line(0,1){8}}

\put(10,9){\oval(10,8)}

\put(12,9){\oval(6,4)}

\put(15,9){\oval(12,12)}

\put(2,15.5){$S$}
\put(3.2,13){$H_2'$}
\put(5.2,10){\small $S \cap C$}
\put(12,9.5){$C'$}
\put(12,11.5){$C$}
\put(17,13.5){$S'$}
\put(.2,16.7){$R$}

\end{picture}
\end{center}
\caption{$k_2$-subspaces of $R$.}\label{figure-subspaces}
\end{figure}
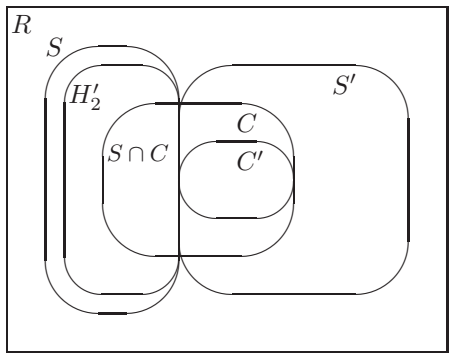

To extend $H_2'$, consider $S, C, H_2',$ and $R$ just as $k_2$-vector spaces as in Fig. ~\ref{figure-subspaces}.  As in the picture below, decompose $C$ as a $k_2$-vector space into \(C = (S \cap C) \oplus C'\).  Then since \(C' \cap S = 0\), we can extend $C'$ to a $k_2$-vector space \(S' \supseteq C'\) such that \(R = S \oplus S'\).  Define \(H_2 := H_2' \oplus S'\), a $k_2$-hyperplane of $R$ since \(\dim_{k_2}(S/H_2') = 1\).  Moreover, $H_2$ contains no nonzero left ideals, since any nonzero left ideal \({_R} L \subseteq H_2\) would contain a minimal left ideal \({_R} L' \subseteq H_2 \cap S = H_2'\).  Most importantly, $H_2$ contains $C$.  

We now define a $k_2$-functional \(\lambda_2: R \rightarrow k_2\) with \(\ker \lambda_2 = H_2\).  Then the \nakaut\ $\sigma_2$ of $R$ as a $k_2$-algebra is defined (\(\forall r \in R\)) by
\[
rs - s \sigma_2(r) \in \ker \lambda_2 = H_2 \all s \in R).  
\]
But since \(rs - s\sigma_1(r) \in C \subseteq H_2\), we have \(\sigma_2(r) = \sigma_1(r)\) for all $r \in R$.  This concludes the proof of Part II\@.
  \qed

\section{Corollaries}

We can now answer the questions posed in the introduction.  We begin with a theorem that does not require the Frobenius assumption.  Let $R$ be a ring that is a finite-dimensional algebra over two fields $k_1$ and $k_2$.  We denote by \lcat\ and \rcat\ the categories of \leftmods\ and \rightmods\ respectively.

Let $\F_i: \lcat \rightarrow \rcat$ be the $k_i$-dual functor:  \(\F_i(\R X) = \rXdi := \Hom{k_i}(X, k_i)\).  Let \biRdi\ be the bimodule \(\Hom{k_i}(R, k_i)\).

\begin{theorem}\label{theorem-bifun}
\(\biRda \simeq \biRdb\) as bimodules iff the functors $\F_1$ and $\F_2$ are naturally equivalent.  
\end{theorem}

\textit{Proof}.
By Brauer's Equivalence Theorem (16.70 in \lomar), the functor $\F_i$ is naturally equivalent to the functor \(\G_i:= \Hom{R}(-,\lRdi)\) on \leftmods, proving the forward direction.  The converse is essentially identical to Theorem 16.71 in \lomar.  We apply the equivalence \(\G_1 \simeq \G_2\)
to the \leftmod\ homomorphism \(\rho_r: \R R \rightarrow \R R\), where \(\rho_r\) is right multiplication by some fixed $r \in R$, as in Fig. ~\ref{figure-brauer}.  Then the map 
\[
\G_i(\rho_r) : \Hom{R}(\R R,\lRdi) \rightarrow \Hom{R}(\R R,\lRdi) 
\]
takes $\alpha$ to the map \( (s \mapsto \alpha (sr)) = r\alpha\), so \(\G_i (\rho_r)\) is {\em left} multiplication by $r$ on \(\Hom{R}(\R R,\lRdi)\).  This gives us a commutative diagram of \rightmods\ as in Fig.~\ref{figure-brauer}.

\setlength{\unitlength}{1in}
\begin{figure}

\begin{center}
\begin{picture}(4.0,1.0)(0,0)

\put(0,0){\(\G_1(\R R)=\Hom{R}(\R R,\lRda)\)}

\put(0,.85){\(\G_1(\R R)=\Hom{R}(\R R,\lRda)\)}

\put(2.2,0){\(\G_2(\R R)=\Hom{R}(\R R,\lRdb)\)}

\put(2.2,.85){\(\G_2(\R R)=\Hom{R}(\R R,\lRdb)\)}

\put(1.8,.05){\vector(1,0){.3}}
\put(1.9,.05){$\sim$}

\put(1.8,.9){\vector(1,0){.3}}
\put(1.9,.9){$\sim$}

\put(.6,.75){\vector(0,-1){.5}}
\put(.7,.45){\(\G_1(\rho_r) = r \cdot\)}
\put(2.8,.75){\vector(0,-1){.5}}
\put(2.9,.45){\(\G_2(\rho_r) = r \cdot\)}

\end{picture}
\end{center}
\caption{The equivalence \(\G_1 \simeq \G_2\) applied to \(\rho_r: \R R \rightarrow \R R\).}\label{figure-brauer}
\end{figure}

However, \(\Hom{R}(\R R,\lRdi) \simeq \rRdi\) as \rightmods\ under the isomorphism \(\alpha \mapsto \alpha(1) \), so the isomorphism on the top and bottom rows is \(\rRda \simeq \rRdb\).  The commutativity of the diagram shows that this isomorphism respects the left $R$-action as well, so we have \(\biRda \simeq \biRdb\) as bimodules, as desired.  
  \qed

To apply this theorem, let \sig\ be any automorphism of $R$ and let $M_R$ be a \rightmod.  We define the twisted \rightmod\ $M_{R^\sigma}$ to be the same abelian group as $M$ with the $R$-action defined by
\[
m*r := m\sigma(r) \all r \in R, m \in M).
\]

(Thanks to Mark Davis for suggesting this 
definition.)  Now let $\R X$ be a \leftmod\ with $k$-dual \(\rXd := \Hom{k} (X,k)\).  Let $(X^*)_R$ denote the $R$-dual $\Hom{R}(\R X,\R R),$ the isomorphism type of which is, of course, independent of $k$.

\begin{theorem}\label{theorem-twistmodule}  Let $R$ be a Frobenius $k$-algebra with Nakayama automorphism \sig.  Then there is a natural \rightmod\ isomorphism \(\rXd \isom (X^*)_{R^\sigma}\).
\end{theorem}

\textit{Proof}.
We have an isomorphism \(\R R \simeq \lRd\), say given by \(1 \mapsto \lambda\).  
Then \(\forall r \in R\),

\[
(\lambda r)(s) = \lambda(rs) = \lambda(s \sigma(r)) = (\sigma(r)\lambda)(s)\all s \in R), 
\]

so \(\lambda r = \sigma(r) \lambda\) in \Rd.  Now By Brauer's Theorem, we have a natural isomorphism \(\rXd \simeq \Hom{R}(\R X, \lRd) \) of \rightmods.  
The isomorphism \(\R R \simeq \lRd\) of \leftmods\ then gives us an abelian group isomorphism
\(
\Hom{R}(\R X, \R R) \simeq \Hom{R}(\R X, \lRd),
\)
which we denote by \(\alpha \mapsto \widehat{\alpha}\).  Then \(\widehat{\alpha}\) is given by 
\begin{equation}\label{alphahat}
\widehat{\alpha}(x) =  (\alpha(x))\lambda \in \Hat{R}.
\end{equation} 

We claim that although `` $\widehat{}$ '' is not in general an isomorphism of \rightmods,  it satisfies \(\widehat{\alpha r} = \widehat{\alpha}\sigma\inv(r)\).  The theorem then follows by identifying \Xd\ with \(\Hom{R}(\R X, \lRd)\) and taking \(f: \Xd \rightarrow \Hom{R}(\R X, \R R)\) to be the inverse of `` $\widehat{}$ ''.  

To prove the claim, let \(x \in X, r \in R\).  Then in \Rd, we have
\begin{eqnarray*}
\widehat{\alpha r}(x) & = &  ((\alpha r) (x)) \lambda \mbox{ by Eq. \ref{alphahat}} \\
& = & (\alpha(x) r) \lambda  \mbox{ by the $R$-action on } \Hom{R}(\R X, \R R) \\
& = & \alpha (x) (r \lambda) \mbox{ by the associativity of the $R$-action on } \lRd \\
& = & \alpha (x) (\lambda \sigma\inv(r)) \mbox{ as shown above} \\
& = & (\alpha (x) \lambda) \sigma\inv(r) \mbox{ by associativity again} \\
& = & (\widehat{\alpha}(x)) \sigma\inv(r) \mbox{ by Eq. \ref{alphahat}} \\
& = & (\widehat{\alpha} \sigma\inv(r))(x) \mbox{ by the $R$-action on }\Hom{R}(\R X, \lRd).
\end{eqnarray*}
So 
\(
\widehat{\alpha r} = \widehat{\alpha}\sigma\inv(r),
\)
proving our claim and the theorem.
  \qed

\begin{corollary}\label{corollary-funind}
If $R$ is a Frobenius $k$-algebra, then the $k$-dual functor \(\F:=\Hom{k}(-,k): \lcat \rightarrow \rcat\) is independent of $k$.
\end{corollary}

\textit{Proof}.
Apply Theorems~\ref{theorem-nakind} and ~\ref{theorem-twistmodule}.
  \qed

\begin{corollary}\label{corollary-biind}
If $R$ is a Frobenius $k$-algebra, then the bimodule isomorphism type of \biRd\ is independent of $k$.  
\end{corollary}

\textit{Proof}.
Apply Corollary~\ref{corollary-funind} and Theorem~\ref{theorem-bifun}.
  \qed

Corollary~\ref{corollary-biind} suggests that there should be a ring-theoretic characterization of \biRd\ as a bimodule analogous to the fact that \(\rRd \isom E((\Rrad)_R)\) as \rightmods.  We do not yet have such a characterization.  

\begin{corollary}\label{corollary-symind}
If $R$ is any finite-dimensional $k$-algebra, then the property of $R$ being a symmetric $k$-algebra is independent of $k$.
\end{corollary}

\textit{Proof}.
We have seen that the question of whether $R$ is a Frobenius $k$-algebra is independent of $k$.  Now apply Corollary~\ref{corollary-biind}.
  \qed

\section {Ring-theoretic characterization of symmetric algebras}

We have seen in Corollary~\ref{corollary-symind} that the property of a $k$-algebra being symmetric is independent of $k$, suggesting that it should be equivalent to a ring-theoretic property.  In the local case, we saw in Corollary~\ref{corollary-localsym} that an algebra is symmetric iff its left socle is not contained in the commutators.  In the general case, we have ring-theoretic conditions for symmetry if we assume that the ground field $k$ is infinite.  

We continue to assume that $R$ is a finite-dimensional algebra over a field $k$.  As before, let \(J = \mbox{rad } R\) be the Jacobson radical and $\Rbar = R/J$.  The following theorem is similar to Theorem 16.14 in \lomar, which states that $R$ is Frobenius iff \(\rsoc \isom \Rbar_R\) and \(\lsoc \isom \R \Rbar\).  We use $S$ to denote $\lsoc$.  

\begin{theorem}\label{theorem-symring}
Suppose $k$ is infinite.  Then $R$ is a symmetric $k$-algebra iff \(\R S_R \simeq \R \Rbar_R\) as $(R,R)$-bimodules and $[R,R]$ contains no nonzero left ideals of $R$.
\end{theorem}

\textit{Proof}.
If $R$ is symmetric, then we have a bimodule isomorphism \(\hi: \R R_R \isommap \biRd\).  Considering $\hi$ as an isomorphism of \leftmods\ and restricting it to $S$, we have an isomorphism \(\hi: S \isommap \mbox{soc}(\lRd)\).  Note, however, that $\hi$ still respects the right action of $R$ on $S$ and $\mbox{soc}(\lRd)$.  By Example 3.41 in \lomar, we have \(\mbox{soc}(\lRd) = \{f \in \Rd: f(J) = 0\}\), which is isomorphic as an $(R,R)$-bimodule to \(\Hom{k}(\Rbar, k)\).  But since \Rbar\ is a semisimple $k$-algebra, hence symmetric, \(\Hom{k}(\Rbar, k) \isom \Rbar\) as $(\Rbar, \Rbar)$-bimodules and hence also as $(R,R)$-bimodules.  Composing all these, we have an $(R,R)$-bimodule isomorphism $S \simeq \Rbar$.

The condition on $[R,R]$ follows from Theorem~\ref{symdef}, which gives us a $k$-linear functional \(\lambda: R \rightarrow k\) such that \([R,R] \subseteq \ker \lambda\), yet $\ker \lambda$ contains no nonzero left ideals of $R$.

Conversely, suppose \(\hi:\R S_R \isommap \R \Rbar_R\) as $(R,R)$-bimodules and $[R,R]$ contains no nonzero left ideals of $R$.  We consider \(\hi ([R,R] \cap S) \subset \Rbar\), which contains no nonzero left ideals of \Rbar\ since $[R,R]$ contains no nonzero left ideals of $R$.  Moreover we claim that \([\Rbar, \Rbar] \subseteq \hi ([R,R] \cap S) \).  Indeed, let \(\bar x, \bar y \in \Rbar\) (where $x,y \in R$), and suppose that \(\bar{y} = \hi(b)\) for some \(b \in S\).  Then using the fact that $\hi$ is a bimodule isomorphism, we have 
\[
\bar x \bar y - \bar y \bar x = \bar x \hi(b) - \hi(b) \bar x = \hi (xb - bx) \in \hi([R,S]) \subseteq \hi ([R,R] \cap S).  
\]
So \(\hi ([R,R] \cap S) \) is a $k$-subspace of \Rbar\ containing no nonzero left ideals and containing the commutators in \Rbar.  By the same argument used in the proof of Theorem~\ref{theorem-nakind}, Part II, we can enlarge \(\hi ([R,R] \cap S) \) to a $k$-hyperplane $U' \subset \Rbar$ containing no nonzero left ideals.  (Here we use the fact that $k$ is infinite.)  Then we can pull back to \(H':= \hi^{-1} (U')\), a $k$-hyperplane of $S$ containing \([R,R] \cap S\) but containing no nonzero left ideals of $R$.  Then, again by the same argument used in Theorem~\ref{theorem-nakind} (using $[R,R]$ in place of the $C$ that was used there), we can extend $H'$ to $H$, a $k$-hyperplane of $R$ containing $[R,R]$ but containing no nonzero left ideals.  Then by Theorem~\ref{symdef}, $R$ is a symmetric algebra.
  \qed

I do not know if Theorem~\ref{theorem-symring} holds without the assumption that $k$ is infinite.  The proof of the forward implication did not use this assumption, so that half certainly remains true.  Conversely, an old result by Nakayama (\cite{nak4}) states that for a finite-dimensional algebra $R$ over a field, \(\lsoc \isom \R\Rbar\) iff \(\rsoc \isom \Rbar_R\).  So if \(\R S_R \simeq \R \Rbar_R\) as $(R,R)$-bimodules, then $R$ is certainly Frobenius, but it does not seem obvious whether $R$ must be symmetric.

\small
\noindent {\bf Acknowledgements}  I would like to thank T.Y.\ Lam, Greg Marks, and Florence Newberger for valuable conversations about this work.

\vfill\eject

\normalsize


\begin{thebibliography}1

\bibitem{frobenius}  F.G.\ Frobenius.  Theorie der hyperkomplexen {G}r\"{o}ssen.  In J.-P.\ Serre (Ed.), Gesammelte Abhandlungen, Springer-Verlag, Berlin, 1968.

\bibitem{lomar}  T.\ Y.\ Lam, A First Course in Noncommutative Rings, in: Grad. Texts in Math., Vol. 131, Springer-Verlag, Berlin–Heidelberg–New York, 1991.

\bibitem{lomar}  T.\ Y.\ Lam, Lectures on Modules and Rings, in: Grad. Texts in Math., Vol. 189, Springer-Verlag, Berlin–
Heidelberg–New York, 1999.
	
\bibitem{nak1}  T.\ Nakayama, On Frobeniusean algebras. I, Ann. of Math. 40 (1939) 611–633.

\bibitem{nak2}  T.\ Nakayama, On Frobeniusean algebras. II, Ann. of Math. (2) 42 (1941) 1–21.

\bibitem{nak3}  T.\ Nakayama, On Frobeniusean algebras. III, Japan. J. Math. 18 (1942) 49–65.

\bibitem{nak4}  T.\ Nakayama, Supplementary remarks on Frobeniusean algebras. I, Proc. Japan Acad. 25 (7) (1949) 45–50.

\end{thebibliography}
\end{document}